\documentclass{ps}
\usepackage{amscd,amssymb,graphics}

\usepackage{hyperref}

\usepackage{bbm}

\newcommand{\R}{{\mathbb R}}
\newcommand{\I}{{\mathbb I}}
\newcommand{\C}{{\mathbb C}}

\newcommand{\ve}{{\varepsilon}}
\def\norm #1{{\left\Vert\,#1\,\right\Vert}}
\newcommand{\abs}[1]{\lvert#1\rvert}
\newcommand{\e}{{\epsilon}}
\newcommand{\N}{{\mathbb N}}

\newcommand{\E}{{\mathbb E}}

\begin{document}

\title{Universal consistency of the $k$-NN rule in metric spaces and Nagata dimension}
\thanks{S.K. would like to thank JICA-FRIENDSHIP scholarship (fellowship D-1590283/ J-1593019) for supporting her doctoral study at Kyoto University and Department of Mathematics, Kyoto University for supporting the Brazil research trip under the KTGU project.}
\thanks{V.G.P. wants to acknowledge support from CNPq (bolsa Pesquisador Visitante, processo 310012/2016) and CAPES (bolsa Professor Visitante Estrangeiro S\^enior, processo 88881.117018/2016-01).}
\runningtitle{$k$-NN rule in metric spaces}
\author{Beno\^{i}t Collins}\address{Department of Mathematics, Kyoto University,
       Kitashirakawa Oiwake-cho, Sakyo-ku, 606-8502, Japan {\tt collins@math.kyoto-u.ac.jp}}
\author{Sushma Kumari}\address{Department of Mathematical Engineering, Musashino University,
1 Chome-1-20 Shinmachi, Nishitokyo, Tokyo 202-8585, Japan \\ {\tt  kumari@musashino-u.ac.jp}}
\author{Vladimir G. Pestov}\address{Instituto de Matem\'atica e Estat\'\i stica, Universidade Federal da Bahia,
Ondina, Salvador--BA, 40.170-115, Brazil}
\secondaddress{Department of Mathematics and Statistics,
       University of Ottawa,
       Ottawa, ON K1N 6N5, Canada {\tt vpest283@uOttawa.ca}}
\date{Version as of 12:57 BST, June 14, 2020}
\begin{abstract} 
The $k$ nearest neighbour learning rule (under the uniform distance tie breaking) is universally consistent in every
metric space $X$ that is sigma-finite dimensional in the sense of Nagata. This was pointed out by C\'erou and Guyader (2006) as a consequence of the main result by those authors, combined with a theorem in real analysis sketched by D. Preiss (1971) (and elaborated in detail by Assouad and Quentin de Gromard (2006)). We show that it is possible to 
give a direct proof along the same lines as the original theorem of Charles J. Stone (1977) about the universal consistency of the $k$-NN classifier in the finite dimensional Euclidean space. The generalization is non-trivial because of the distance ties being more prevalent in the non-Euclidean setting, and on the way we investigate the relevant geometric properties of the metrics and the limitations of the Stone argument, by constructing various examples.
 \end{abstract}
\begin{resume} 
La r\`egle d'apprentissage des $k$ plus proches voisins (sous le bris uniforme d'\'egalit\'e des distances) est universellment consistente dans chaque espace m\'etrique s\'eparable de dimension sigma-finie au sens de Nagata. Comme indiqu\'e par C\'erou et Guyader (2006),
le r\'esultat fait suite \`a une combinaison du th\'eor\`eme principal de ces auteurs avec un th\'eor\`eme d'analyse r\'eelle esquiss\'e par D. Preiss (1971) (et \'elabor\'e en d\'etail par Assouad et Quentin de Gromard (2006)). Nous montrons qu'il est possible de
donner une preuve directe dans le m\^eme esprit que le th\'eor\`eme original de Charles J. Stone (1977) sur la consistence universelle du classificateur $k$-NN dans l'espace euclidien de dimension finie. La g\'en\'eralisation est non-triviale, car l'\'egalit\'e des distances est plus commune dans le cas non-euclidien, et pendant l'\'elaboration de notre preuve, nous \'etudions des propri\'et\'es g\'eom\'etriques pertinentes des m\'etriques et testons des limites de l'argument de Stone, en construisant quelques exemples.
\end{resume}
\maketitle
\section*{Introduction}
The $k$-nearest neighbour classifier, in spite of being arguably the oldest supervised learning algorithm in existence, still retains his importance, both practical and theoretical. In particular, it was the first classification learning rule whose (weak) universal consistency (in the finite-dimensional Euclidean space $\R^d=\ell^2(d)$) was established, by Charles J. Stone in \cite{stone:77}.

Stone's result is easily extended to all finite-dimensional normed spaces, see, e.g., \cite{duan}. However, the $k$-NN classifier is no longer universally consistent already in the infinite-dimensional Hilbert space $\ell^2$. A series of examples of this kind, obtained in the setting of real analysis, belongs to Preiss, and the first of them \cite{Preiss_1979} is so simple that it can be described in a few lines. We will reproduce it in the article, since the example remains virtually unknown in the statistical machine learning community.

There is sufficient empirical evidence to support the view that the performance of the $k$-NN classifier greatly depends on the chosen metric on the domain (see e.g., \cite{hatko}). There is a supervised learning algorithm, Large Margin Nearest Neighbour Classifier (LMNN), based on the idea of optimizing the $k$-NN performance over all Euclidean metrics on a finite dimensional vector space \cite{WS}. At the same time, it appears that a theoretical foundation for such an optimization over a set of distances is still lacking.
The first question to address in this connection, is of course to characterize those metrics (generating the original Borel structure of the domain) for which the $k$-NN classifier is (weakly) universally consistent. 

While the problem in this generality remains still open, a great advance in this direction was made by C\'erou and Guyader in \cite{CG}. They have shown that the $k$-NN classifier is consistent under the assumption that the regression function $\eta(x)$ satisfies the weak Lebesgue--Besicovitch differentiation property:
\begin{equation}
\frac{1}{\mu(B_r(x))}\int_{B_r(x)} \eta(x)\,d\mu(x) \to \eta(x),
\label{eq:lb}
\end{equation}
where the convergence is in measure, that is, for each $\e>0$,
\[\mu\left\{x\in\Omega\colon 
\left\vert\frac{1}{\mu(B_r(x))}\int_{B_r(x)} \eta(x)\,d\mu(x) - \eta(x)
\right\vert >\e\right\}\to 0\mbox{ when } r\downarrow 0.\]
The probability measure $\mu$ above is the sample distribution law.
The proof extended the ideas of the paper \cite{devroye}, in which it was previously observed that Stone's universal consistency can be deduced from the classical Lebesgue--Besicovitch differentiation theorem: every $L^1(\mu)$-function $f$ on $\R^d$ satisfies Eq. (\ref{eq:lb}), even in the strong sense (convergence almost everywhere). See also \cite{FFL}.

Those separable metric spaces in which the weak Lebesgue--Besicovitch differentiation property holds for every Borel probability measure (equivalently, for every sigma-finite locally finite Borel measure) have not yet been characterized. 
But the complete separable metric spaces in which the {\em strong} Lebesgue--Besicovitch differentiation property holds for every such measure as above have been described by Preiss \cite{preiss83}: they are exactly those spaces that are {\em sigma-finite dimensional} in the sense of Nagata \cite{nagata,ostrand}. (For finite dimensional spaces in the sense of Nagata, the sketch of a proof by Preiss, in the sufficiency direction, was elaborated by Assouad and Quentin de Gromard in \cite{AQ}. The completeness assumption on the metric space is only essential for the necessity part of the result.) In particular, it follows that every sigma-finite dimensional separable metric space satisfies the weak Lebesgue--Besicovitch differentiation property for every probability measure.

Combining the result of Preiss with that of C\'erou--Guyader, one concludes that the $k$-NN classifier is universally consistent in every sigma-finite dimensional separable metric space, as was noted in \cite{CG}.

The authors of \cite{CG} mention in their paper that {\em ``[Stone's theorem] is based on a geometrical result, known as Stone's Lemma. This powerful and elegant argument can unfortunately not be generalized to infinite dimension.''} The aim of this article is to show that at least Stone's original proof, including Stone's geometric lemma as its main tool, can be extended from the Euclidean case to the sigma-finite dimensional metric spaces. In fact, as we will show, the geometry behind Stone's lemma, even if it appears to be essentially based on the Euclidean structure of the space, is captured by the notion of Nagata dimension, which is a purely metric concept. In this way, Stone's geometric lemma and indeed the original Stone's proof of the universal consistency of the $k$-NN classifier, become applicable to a wide range of metric spaces.

In the absence of distance ties (that is, in case where every sphere is a $\mu$-negligible set with regard to the underlying measure $\mu$), the extension is quite straightforward, indeed almost literal. However, this is not so in the presence of distance ties: an example shows that the conclusion of Stone's geometric lemma may not hold. Another example shows that even in the compact metric spaces of Nagata dimension zero, the distance ties may be everpresent. We also show that an attempt to reduce the case to the situation without distance ties by learning in the product of $\Omega$ with the unit interval (an additional random variable used for tie-breaking) cannot work, because already the product of a zero-dimensional space in the sense of Nagata with the interval (which has dimension one) can have an infinite Nagata dimension. Stone's geometric lemma has to be modified, to parallel the Hardy--Littlewood inequality in the geometric measure theory. 

We do not touch upon the subject of strong universal consistency in general metric spaces. 
The main open question left is whether every metric space in which the $k$-NN classifier is universally consistent is necessarily sigma-finite dimensional. A positive answer, modulo the work of \cite{CG} and \cite{preiss83}, would also answer in the affirmative an open question in real analysis going back to Preiss: suppose a metric space $X$ satisfies the weak Lebesgue--Besicovitch differentiation property for every sigma-finite locally finite Borel measure, will it satisfy the strong Lebesgue--Besicovitch differentiation property for every such measure?

\section{Setting for statistical learning}
Here we will recall the standard probabilistic model for statistical learning theory.
The {\em domain}, $\Omega$, means a standard Borel space, that is, a set equipped with a sigma-algebra which coincides with the sigma-algebra of Borel sets generated by a suitable separable complete metric. (Recall that the Borel structure generated by a metric $\rho$ on a set $\Omega$ is the smallest sigma-algebra containing all open subsets of the metric space $(\Omega,\rho)$.) The distribution laws for datapoints, both unlabelled and labelled, are Borel probability measures defined on the corresponding Borel sigma-algebra.

Since we will be dealing with the $k$-NN classifier, the domain, $\Omega$, will actually be a metric space, which we also assume to be separable.

Labelled data pairs $(x,y)$, where $x\in\Omega$ and $y\in\{0,1\}$, will follow an unknown probability distribution $\tilde\mu$, that is, a Borel probability measure on $\Omega\times\{0,1\}$. We denote the corresponding random element $(X,Y)\sim \tilde\mu$. Define two Borel measures on $\Omega$, $\mu_i$, $i=0,1$, by $\mu_i(A)=\tilde\mu(A\times\{i\})$. In this way, $\mu_0$ is governing the distribution of the elements labelled $0$, and similarly for $\mu_1$. The sum $\mu=\mu_0+\mu_1$ (the direct image of $\tilde\mu$ under the projection from $\Omega\times\{0,1\}$ onto $\Omega$) is a Borel probability measure on $\Omega$, the distribution law of unlabelled data points. Clearly, $\mu_i$ is absolutely continuous with regard to $\mu$, that is, if $\mu(A)=0$, then $\mu_i(A)=0$ for $i=0,1$. The corresponding Radon-Nikod\'ym derivative in the case $i=1$ is just the conditional probability for a point $x$ to be labeled $1$:
\[\eta(x) = \frac{d\mu_1}{d\tilde\mu}(x) = P[Y=1\vert X=x].\]
In statistical terminology, $\eta$ is the {\em regression function.}

\begin{figure}
   \begin{center}
    \scalebox{0.25}{\includegraphics{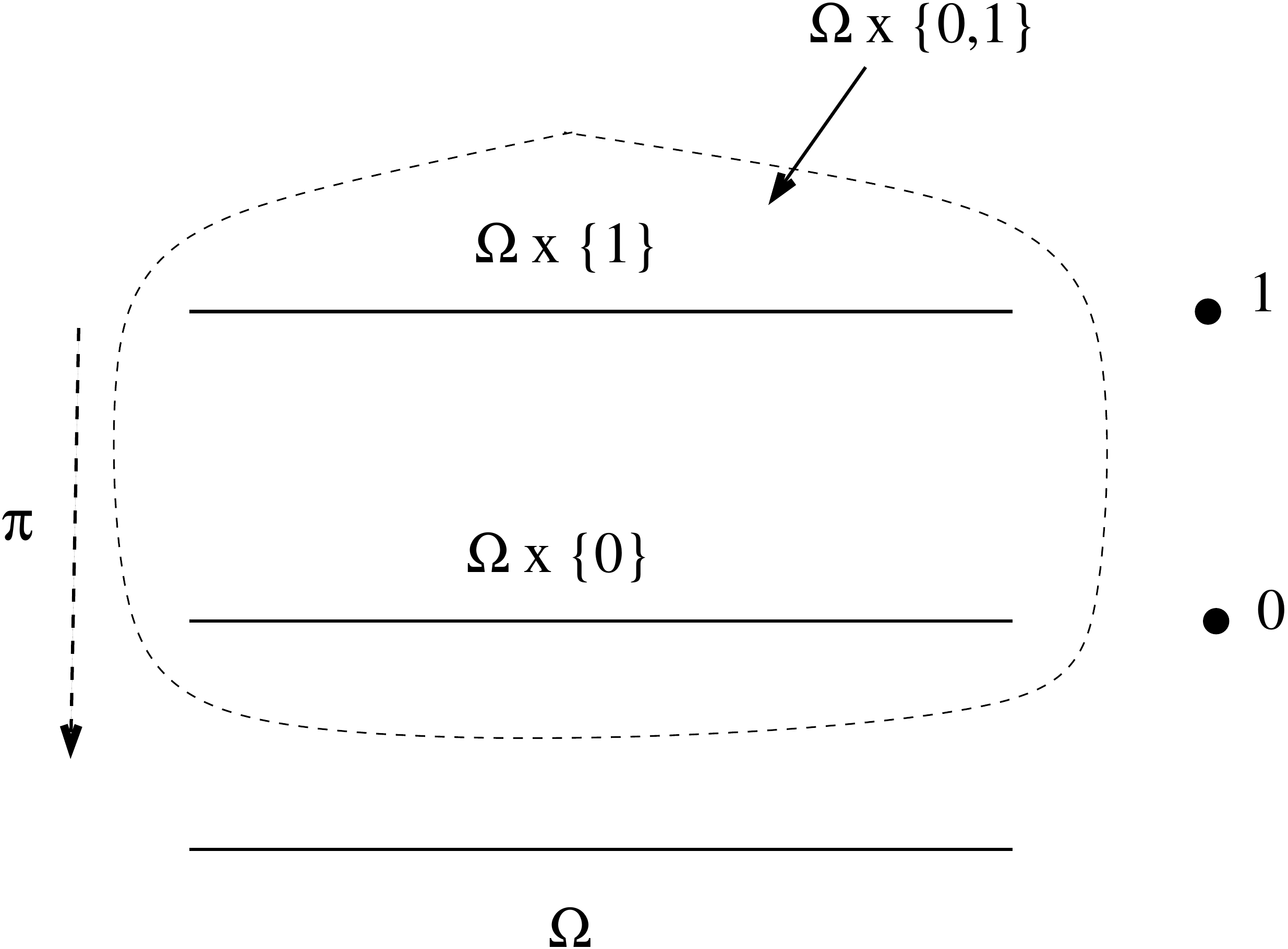}} 
  \end{center}
\caption{Labeled domain and the projection $\pi\colon\Omega\times\{0,1\}\to\Omega$.}
\end{figure}

Together with the Borel probability measure $\mu$ on $\Omega$, the regression function allows for an alternative, and often more convenient, description of the joint law $\tilde\mu$. Namely, given $A\subseteq\Omega$,
\[\mu_1(A)=\int_A \eta(x)\,d\mu,\]
and
\[\mu_0(A)=\int_A (1-\eta(x))\,d\mu,\]
which allows to reconstruct the measure $\tilde\mu$ on $\Omega\times\{0,1\}$.

Let ${\mathcal B}(\Omega,\{0,1\})$ denote the collection of all Borel measurable binary functions on the domain, that is, essentially, the family of all Borel subsets of $\Omega$.
Given such a function $f\colon\Omega\to\{0,1\}$ (a {\em classifier}), the  
{\em misclassification error} is defined by 
\[{\mathrm{err}}_{\tilde\mu}(f) = \tilde\mu\{(x,y)\colon f(x) \neq y\}=P[f(X)\neq Y].\]
The {\em Bayes error} is the infimal misclassification error taken over all possible classifiers: 
\[\ell^{\ast}=\ell^{\ast}(\tilde\mu)=\inf_{f}{\mathrm{err}}_{\tilde\mu}(f).\]
It is a simple exercise to verify that the Bayes error is achieved on some classifier (and thus is the minimum), which is called a {\em Bayes classifier}. For instance, every classifier satisfying
\[T_{bayes}(x) = \begin{cases} 1,&\mbox{ if }\eta(x)>\frac 12,\\
0,&\mbox{ if }\eta(x)<\frac 12,
\end{cases}\]
is a Bayes classifier. 

The Bayes error is zero if and only if the learning problem is deterministic, that is, the regression function $\eta$ is equal almost everywhere to the indicator function, $\chi_C$, of a {\em concept} $C\subseteq\Omega$, a Borel subset of the domain.

A {\em learning rule} is a family of mappings ${\mathcal L}=\left({\mathcal L}_n\right)_{n=1}^{\infty}$, where
\[{\mathcal L}_n\colon \Omega^n\times \{0,1\}^n \to {\mathcal B}(\Omega,\{0,1\}),~~n=1,2,\ldots
\]
and the functions ${\mathcal L}_n$ satisfy the following measurability assumption: the associated maps
\[\Omega^n\times \{0,1\}^n\times \Omega \ni (\sigma,x)\mapsto {\mathcal L}_n(\sigma)(x) \in \{0,1\}\]
are Borel (or just universally measurable). Here $\sigma=(x_1,\ldots,x_n,y_1,\ldots,y_n)$ is a labelled learning sample.

The data is modelled by a sequence of independent identically distributed random elements $(X_n,Y_n)$ of $\Omega\times\{0,1\}$, following the law $\tilde\mu$. 
Denote $\varsigma$ an infinite sample path. In this context, $\mathcal L_n$ only gets to see the first $n$ labelled coordinates of $\varsigma$.
A learning rule $\mathcal L$ is {\em weakly consistent}, or simply {\em consistent,} if ${\mathrm{err}}_{\mu}{\mathcal L}_n(\varsigma)\to \ell^{\ast}$ in probability as $n\to\infty$. If the convergence occurs almost surely (that is, along almost all sample paths $\varsigma\sim \tilde\mu^{\infty}$), then $\mathcal L$ is said to be {\em strongly consistent}.
Finally, $\mathcal L$ is {\em universally (weakly / strongly) consistent} if it is weakly / strongly consistent under every Borel probability measure $\tilde\mu$ on the standard Borel space $\Omega\times\{0,1\}$. 

The learning rule we study is the $k$-NN classifier, defined by selecting the label 
${\mathcal L}_n(\sigma)(x)\in\{0,1\}$ by the majority vote among the values of $y$ corresponding to the $k=k_n$ nearest neighbours of $x$ in the learning sample $\sigma$. 

If $k$ is even, then a voting tie may occur. This is of lesser importance, and can be broken in any way. For instance, by always assigning the value $1$ in case of a voting tie, or by choosing the value randomly. The consistency results usually do not depend on it. Intuitively, if voting ties keep occurring asymptotically at a point $x$ along a sample path, it means that $\eta(x)=1/2$ and so any value of the classifier assigned to $x$ would do.

It may also happen that the smallest closed ball containing $k$ nearest neighbours of a point $x$ contains more than $k$ elements of a sample (distance ties). This situation is more difficult to manage and requires a consistent tie-breaking strategy, whose choice may affect the consistency results. 

Given $k$ and $n\geq k$, we define $r^{\varsigma_n}_{k\mbox{\tiny -NN}}(x)$ as the smallest radius of a closed ball around $x$ containing at least $k$ nearest neighbours of $x$ in the sample $\varsigma_n$:
\begin{equation}
r^{\varsigma_n}_{k\mbox{\tiny -NN}}(x)=\min\{r\geq 0\colon \sharp\{i=1,2,\ldots,n\colon x_i\in \bar B_r(x)\}\geq k\}.
\label{eq:rknn}
\end{equation}
As the corresponding open ball around $x$ contains at most $k-1$ elements of the sample, the ties may only occur on the sphere.

We adopt the combinatorial notation $[n]=\{1,2,\ldots n\}$.
If $\sigma\in\Omega^n$ and $\sigma^\prime\in\Omega^k$, $k\leq n$, the symbol
\[\sigma^\prime\sqsubset\sigma\]
means that there is an injection $f\colon [k]\to [n]$ such that
\[\forall i=1,2,\ldots,k,~\sigma^\prime_i=\sigma_{f(i)}.\] 
A {\em $k$ nearest neighbour map} is a function 
\[k\mbox{-NN}^{\sigma}\colon\Omega^n\times\Omega\to\Omega^k\]
with the properties 
\begin{enumerate}
\item $k\mbox{-NN}^{\sigma}(x) \sqsubset \sigma$, and
\item all points $x_i$ in $\sigma$ that are at a distance strictly less than $r^{\varsigma_n}_{k\mbox{\tiny -NN}}(x)$ to $x$ are in $k\mbox{-NN}^{\sigma}(x)$. 
\end{enumerate}

The mapping $k\mbox{-NN}^{\sigma}$ can be deterministic or stochastic, in which case it will depend on an additional random variable, independent of the sample.

An example of the former kind is based on the natural order on the sample, $x_1<x_2<\ldots<x_n$. In this case, from among the points belonging to the sphere of radius $r_{k\mbox{\tiny -NN}^{\sigma}}(x)$ around $x$ we choose points with the smaller index: $k\mbox{-NN}^{\sigma}(x)$ contains all the points of $\sigma$ in the open ball, $B_{r_{k\mbox{\tiny -NN}^{\sigma}}(x)}(x)$, plus a necessary number (at least one) of points of $\sigma\cap S_{r_{k\mbox{\tiny -NN}^{\sigma}}(x)}(x)$ having smallest indices.

An example of the second kind is to use a similar procedure, after applying a random permutation of the indices first. A random learning input will consist of a pair $(W_n,P_n)$, where $W_n$ is a random $n$-sample and $P_n$ is a random element of the group of permutations of rank $n$. An equivalent (and more common) way would be to use a sequence of i.i.d. random elements $Z_n$ of the unit interval or the real line, distributed according to the uniform (resp. gaussian) law, and in case of a tie, give a preference to a realization $x_i$ over $x_j$ provided the value $z_i$ is smaller than $z_j$.

Now, a formal definition of the $k$-NN learning rule can be given as follows:
\begin{eqnarray*}
{\mathcal L}^{k\mbox{\tiny -NN}}_n(\sigma,\e)(x) &=&
\chi_{[0,\infty)}\left[\frac 1k\sum_{x_i\in k\mbox{\tiny -NN}^{\sigma}(x)}\e_i-\frac 12 \right]
\\[1mm]
&=& 
\chi_{[0,\infty)}\left[\E_{\mu_{k\mbox{\tiny -NN}^{\sigma}(x)}}\e -\frac 12\right].
\end{eqnarray*}
Here, $\chi_{[0,\infty)}$ is the Heaviside function, the sign of the argument:
\[\chi_{[0,\infty)}(t) =\begin{cases} 1,&\mbox{ if }t\geq 0,\\
0,&\mbox{ if }t<0.\end{cases}\]
The empirical measure $\mu_{k\mbox{\tiny -NN}^{\sigma}(x)}$ is a uniform measure supported on the set of $k$ nearest neighbours of $x$ within the sample $\sigma$, and the label $\e$ is seen as a function $\e\colon \{x_1,x_2,\ldots,x_n\}\to\{0,1\}$.

The expression appearing under the argument,
\begin{equation}
\eta_{n,k} = \frac 1k\sum_{x_i\in k\mbox{\tiny -NN}^{\sigma}(x)}\e_i,
\label{eq:empregfun}
\end{equation}
is the {\em empirical regression function.} In the presence of a law of labelled points, it is a random variable, and so we have the following immediate, yet important, observation.

\begin{prpstn}
Let $(\mu,\eta)$ be a learning problem in a separable metric space $(\Omega,d)$.
If the values of the empirical regression function, $\eta_{n,k}$, converge to $\eta$ in probability (resp. almost surely) in the region
\[\Omega_{\eta} = \left\{x\in\Omega\colon \eta(x)\neq \frac 12\right\},\] 
then the $k$-NN classifier is consistent (resp. strongly consistent) under $(\mu,\eta)$. 
\label{p:empirical_regression}
\end{prpstn}

We conclude this section by recalling an important technical tool.

\begin{thrm}[Cover-Hart lemma \cite{CH}]
Let $\Omega$ be a separable metric space, and let $\mu$ be a Borel probability measure on $\Omega$. Almost surely, the function $r^{\varsigma_n}_{k\mbox{\tiny -NN}}$ (Eq. (\ref{eq:rknn})) converges to zero uniformly over any precompact subset $K\subseteq\mbox{supp}\,\mu$.
\label{l:cover-hart}
\end{thrm}

\begin{proof} 
Let $A$ be a countable dense subset of $\mbox{supp}\,\mu$. A standard argument shows that, almost surely, for all $a\in A$ and each rational $\e>0$, the open ball $B_{\e}(a)$ contains an infinite number of elements of a sample path. Consequently, the functions $r^{\varsigma_n}_{k\mbox{\tiny -NN}}\colon\Omega\to\R$ almost surely converge to zero pointwise on $A$ as $n\to\infty$. Since these functions are easily seen to be $1$-Lipschitz and in particular form a uniformly equicontinuous family, we conclude.
\end{proof}

\section{Example of Preiss\label{s:preiss}}

Here we will discuss a 1979 example of Preiss \cite{Preiss_1979}. Preiss's aim was to prove that the Lebesgue--Besicovitch differentiation property (Eq. (\ref{eq:lb})) fails in the infinite-dimensional Hilbert space $\ell^2$. However, as already suggested in \cite{CG}, his example can be easily adapted to prove that the $k$-NN learning rule is not universally consistent in the infinite-dimensional separable Hilbert space $\ell^2$ either.

Recall the notation $[n]=\{1,2,\ldots n\}$. Let $(N_k)$ be a sequence of positive natural numbers $\geq 2$, to be selected later.
Denote by
\[Q=\prod_{k=1}^{\infty} [N_k]\] 
the Cartesian product of finite discrete spaces equipped with the product topology. It is a Cantor space (the unique, up to a homeomorphism, totally disconnected compact metrizable space without isolated points). 

Let $\pi_k$ denote the canonical cordinate projections of $Q$ on the $k$-dimensional cubes $Q_k=\prod_{i=1}^{k} [N_i]$. 
Denote $Q^\ast= \cup_{k=1}^{\infty} Q_k$ a disjoint union of the cubes $Q_k$, 
and let ${\mathcal H}=\ell^2(Q^\ast)$ be a Hilbert space spanned by an orthonormal basis $(e_{\bar n})$ indexed by elements $\bar n$ of this union.

For every $\bar n=(n_1,\ldots,n_i,\ldots)\in Q$ define
\[f(\bar n) = \sum_{i=1}^{\infty} 2^{-i}e_{(n_1,\ldots,n_i)}\in {\mathcal H}.\]
The map $f$ is continuous and injective, thus a homeomorphism onto its image. Denote $\nu$ the Haar measure on $Q$ (the product of the uniform measures on all $[N_k]$).
Let $\mu_1=f_{\ast}(\nu)$ be the direct image of $\nu$, a compactly-supported Borel probability measure on $\mathcal H$. If $r>0$ satisfies $2^{-k}\leq r^2<2^{-k+1}$, then for each $\bar n =(n_1,n_2,\ldots)\in Q$,
\begin{align*}
\mu_1(B_r(f(\bar n))) & = \nu(\pi_{k+1}^{-1}(\bar n)) \\
& = (N_1N_2\ldots N_{k+1})^{-1}.
\end{align*}

Now, for every $k$ and each $\bar n=(n_1,\ldots,n_k)\in Q_k\subseteq Q^\ast$ define in a similar way
\[f(\bar n)=\sum_{i=1}^k 2^{-i}e_{(n_1,\ldots,n_i)}\in {\mathcal H}.\]
Note that the closure of $f(Q^\ast)$ contains $f(Q)$ (as a proper subset). 
Now define a purely atomic measure $\mu_0$ supported on the image of $Q^\ast$ under $f$, having the following special form:
\[\mu_0=\sum_{k=1}^{\infty}\sum_{\bar n\in Q_k} a_k\delta_{\bar n}.\]
The weights $a_k>0$ are chosen so that the measure is finite:
\begin{equation}
\sum_{k=1}^{\infty} a_k <\infty.
\label{eq:finite}
\end{equation}
Since for $r$ satisfying $2^{-k}\leq r^2<2^{-k+1}$ and $\bar n\in Q$ the ball $B_r(f(\bar n))$ contains in particular $f(n_1n_2,\ldots,n_k)$, we have
\[\mu_0(B_r(f(\bar n))) \geq a_k.\]
Assuming in addition that 
\begin{equation}
a_kN_1N_2\ldots N_kN_{k+1} \to \infty\mbox{ as }k\to\infty,
\label{eq:infinity}
\end{equation}
we conclude:
\[\frac{\mu_1(B_r(f(\bar n)))}{\mu_0(B_r(f(\bar n)))}\leq\frac{(N_1N_2\ldots N_{k+1})^{-1}}{a_k}\to 0\mbox{ when } r\downarrow 0.\]
Clearly, the conditions (\ref{eq:finite}) and (\ref{eq:infinity}) can be simultaneously satisfied by a recursive choice of $(N_k)$ and $(a_k)$. 

Now renormalize the measures $\mu_0$ and $\mu_1$ so that $\mu=\mu_0+\mu_1$ is a probability measure, and interpret $\mu_i$ as the distribution of points labelled $i=0,1$. Thus, the regression function is deterministic, $\eta=\chi_C$, where we are learning the concept $C=f(Q)={\mathrm{supp}}\,\mu_1$, $\mu_1(C)>0$. 

For a random element $X\in {\mathcal H}$, $X\sim\mu$, the distance $r_k(X)$ to the $k$-th nearest neighbour within an i.i.d. $n$-sample goes to zero almost surely when $k/n\to 0$, according to a lemma of Cover and Hart, and the convergence is uniform on the precompact support of $\mu$.
It follows that the probability of one of the $k$ nearest neighbours to a random point $X\in {\mathcal H}$ to be labelled one, conditionally on $r^{\varsigma_n}_{k\mbox{\tiny -NN}}=r$, converges to zero, uniformly in $r$. The $k$-NN learning rule will almost surely predict a sequence of classifiers converging to the identically zero classifier, and so is not consistent.

\section{Classical theorem of Charles J. Stone}

\subsection{The case of continuous regression function}
Proposition \ref{p:empirical_regression} and the Cover--Hart lemma \ref{l:cover-hart} together imply that the $k$-NN classifier is universally consistent in a separable metric space whenever the regression function $\eta$ is continuous. In view of Proposition \ref{p:empirical_regression}, it is enough to make the following observation.

\begin{lmm} 
Let $(\Omega,\mu)$ be a separable metric space equipped with a Borel probability measure, and let $\eta$ be a continuous regression function. Then 
\[\eta_{n,k} \to \eta\]
in probability, when $n,k\to\infty$, $k/n\to 0$.
\label{l:etaempirica}
\end{lmm}

\begin{proof}
It follows from the Cover--Hart lemma that
the set $k\mbox{-NN}^{\varsigma_n}(x)$ of $k$ nearest neighbours of $x$ almost surely converges to $x$, for almost all $x\in\mbox{supp}\,\mu$, and since $\eta$ is continuous, the set of values $\eta(k\mbox{-NN}(x))$ almost surely coverges to $\eta(x)$ in an obvious sense: for every $\ve>0$, there exists $N$ such that
\begin{equation}
\forall n\geq N,~~\eta(k\mbox{-NN}(x))\subseteq (\eta(x)-\ve,\eta(x)+\ve),
\label{eq:veN}
\end{equation}
where $k$ depends on $n$.
Let $\ve>0$ and $N$ be fixed, and denote $P_{\ve,N}$ the set of pairs $(\varsigma,x)$ consisting of a sample path $\varsigma$ and a point $x\in\Omega$ satisfying Eq. (\ref{eq:veN}). Select $N$ with the property $\mu(P_{\ve,N})>1-\ve$.
Let $\e=(\e_i)_{i=1}^{\infty}$ denote the sequence of labels for $\varsigma$, which is a random variable with the joint law $\otimes_{n=1}^{\infty} \{\eta(x_i),1-\eta(x_i)\}$. By the above, whenever $(\varsigma,x)\in P_{\ve,N}$ and $n\geq N$, if $x_i$ is one of the $k$ nearest neighbours of $x$ in $\varsigma_n$, we have $\E\e_i=\eta(x_i)\in (\eta(x)-\ve,\eta(x)+\ve)$.
According to a version of the Law of Large Numbers with Chernoff's bounds, the probability of the event
\begin{equation}
\frac{\sum_{x_i\in k\mbox{\tiny -NN}(x)}\e_i}{k}\notin (\eta(x)-2\ve,\eta(x)+2\ve)
\label{eq:deviationevent}
\end{equation}
is exponentially small, bounded above by $2\exp(-2\ve^2 k)$. Thus, when $n\geq N$, $P[\left\vert\eta_{n,k}-\eta\right\vert\geq\ve]<\ve+2\exp(-2\ve^2 k)$, and we conclude.
\end{proof}

\begin{rmrk}
In the most general case (with the uniform tie-breaking) we can only infer the almost sure convergence if $k=k_n$ grows fast enough as a function in $n$, for otherwise the series $\sum_{n=1}^\infty 2\exp(-2\ve^2 k_n)$ may be divergent. 
\end{rmrk}

\subsection{Stone's geometric lemma for $\R^d$}
In the case of a general Borel regression function $\eta$, which can be discontinuous $\mu$-almost everywhere, where $\mu$, as before, is the sample distribution on a separable metric space, a version of the classical Luzin theorem of real analysis says that for any $\ve>0$ there is a closed precompact set $K$ of measure $\mu(K)>1-\ve$ upon which $\eta$ is continuous. (See Appendix.) Now we have control over the behaviour of those $k$-nearest neighbours of a point $x$ that belong to $K$: the mean value of the regression function $\eta$ taken at those $k$-nearest neighbours will converge to $\eta(x)$. However, we have no control over the behaviour of the values of $\eta$ at the $k$-nearest neighbours of $x$ that belong to the open set $U=\Omega\setminus K$.
The problem is therefore to limit the influence of the remaining $\approx \ve n$ sample points belonging to $U$. Intuitively, as the example of Preiss shows, in infinite dimensions the influence of the few points outside of $K$ can become overwhelming, no matter how close the measure of $K$ is to one.

In the Euclidean case, this goal is achieved with the help of Stone's geometric lemma, which uses the finite-dimensional Euclidean structure of the space in a beautiful way.

\begin{lmm}[Stone's geometric lemma for $\R^d$]
For every natural $d$, there is an absolute constant $C=C(d)$ with the following property. Let
\[\sigma=(x_1,x_2,\ldots,x_n),~x_i\in\R^d,~i=1,2,\ldots,n,\]
be a finite sample in $\ell^2(d)$ (possibly with repetitions), and let $x\in\ell^2(d)$ be any. Given $k\in\N_+$, the number of $i$ such that $x\neq x_i$ and $x$ is among the $k$ nearest neighbours of $x_i$ inside the sample
\begin{equation}
x_1,x_2,\ldots, x_{i-1},x,x_{i+1},\dots,x_n
\label{eq:subst_sample}
\end{equation}
is limited by $Ck$. 
\label{l:stonek=k}
\end{lmm}

\begin{figure}
   \begin{center}
    \scalebox{0.3}{\includegraphics{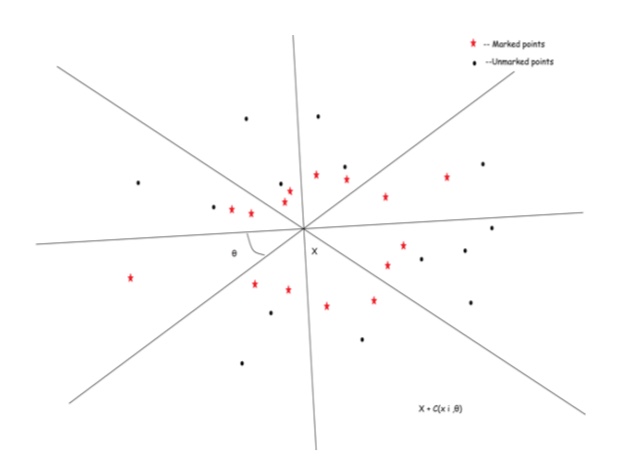}} 
  \end{center}
\caption{To the proof of Stone's geometric lemma (case $k=2$).}
\end{figure}

\begin{proof} 
Cover $\R^d$ with $C=C(d)$ cones of central angle $< \pi/3$ with vertices at $x$. Inside each cone mark the maximal possible number $\leq k$ of the nearest neighbours of $x$, that are set-theoretically different from $x$. (The strategy for possible distance tie-breaking is unimpotant.) In this way, up to $Ck$ points are marked. Let now $i$ be any, such that $x_i\neq x$ as a point. If $x_i$ has not been marked, this means the cone containing $x_i$ has $k$ points, different from $x$, that have been marked. Consider any of the marked points inside the same cone, say $y$. A simple argument of planimetry, inside an affine plane passing through $x$, $x_i$, and $y$, shows that 
\begin{equation}
\norm{x_i-x}>\norm{x_i-y},
\label{eq:stoneequation}
\end{equation}
and so the $k$ nearest neighbours of $x_i$ inside the sample in Eq. (\ref{eq:subst_sample}) will all be among the marked points, excluding $x$.
\end{proof}

\begin{rmrk}
Note that in the statement of Stone's geometric lemma neither the order of the sample $x_1,x_2,\ldots,\allowbreak x_{i-1}, x,x_{i+1},\dots,x_n$, nor the tie-breaking strategy are of any importance.
\label{r:unimportant}
\end{rmrk}

\begin{rmrk}
If the cones have central angle $\pi/3=60^\circ$, then the displayed inequality in the proof of the lemma (Eq. (\ref{eq:stoneequation})) is no longer strict. This is less convenient in case of distance ties.
\label{r:strictineq}
\end{rmrk}

\subsection{Proof of Stone's theorem\label{ss:stone}}

\begin{thrm}[Charles J. Stone, 1977]
Let $k,n\to\infty$, $k/n\to 0$. Then the $k$-NN classification rule in the finite-dimensional Euclidean space $\R^d$ is universally consistent.
\label{t:stone}
\end{thrm}

Let us begin with the case where the domain $\Omega$ is an arbitrary separable metric space, $\mu$ is a probability measure on $\Omega$, and $\eta$ is a regression function. 
Let $\ve>0$ be any. Using a variation of Luzin's theorem (Theorem \ref{th:luzin} in Appendix), choose a closed precompact set $K\subseteq\Omega$ with $\mu(K)>1-\ve$ and $\eta\vert_K$ continuous. Denote $U=\Omega\setminus K$. By the Tietze--Urysohn extension theorem (see e.g., \cite{E}, 2.1.8), the function $\eta\vert_K$ extends to a continuous $[0,1]$-valued function $\psi$ on all of $\Omega$. 

Denote $\psi_{n,k}$ the empirical regression function (Eq. (\ref{eq:empregfun})) corresponding to $\psi$ viewed as a regression function on its own right.
We have:
\begin{align*}
\E\left\vert\eta-\eta_{n,k}\right\vert &\leq
\underbrace{\E\left\vert\eta-\psi\right\vert}_{(I)} + 
\underbrace{\E\left\vert\psi - \psi_{n,k}\right\vert}_{(II)} + \underbrace{\E\left\vert\psi_{n,k}-\eta_{n,k}\right\vert}_{(III)},
\end{align*}
where $(I)\leq \mu(U)<\ve$ and $(II)\to 0$ in probability by virtue of Lemma \ref{l:etaempirica}. It only remains to estimate the term $(III)$. 

With this purpose, let now $\Omega=\ell^2(d)$, let $K$ be a compact subset of $\R^d$, and let $U=\Omega\setminus K$. Let $X,X_1,X_2,\ldots,X_n$ be i.i.d. random elements of $\R^d$ following the law $\mu$.
We will estimate the expected number of the random elements $X_i$ of an $n$-sample that (1) belong to $U$, and (2) are among the $k$ nearest neighbours of a random element $X$ belonging to $K$. Applying the symmetrization with a transposition of the coordinates $\tau_i\colon X\leftrightarrow X_i$, as well as  Stone's geometric lemma \ref{l:stonek=k}, we obtain:

\begin{align*}
\E \frac 1 k \sharp \{i=1,2,\ldots,n\}{\mathbbm{1}}_{X_i\in k\mbox{\tiny -NN}(X)}{\mathbbm{1}}_{X\in K}{\mathbbm{1}}_{X_i\notin K}
& = \E \frac 1 k \sum \chi_U(X_i){\mathbbm{1}}_{X_i\in k\mbox{\tiny -NN}(X)}{\mathbbm{1}}_{X\in K}{\mathbbm{1}}_{X_i\notin K} \\
&\leq \E \frac 1 k \sum \chi_U(X_i){\mathbbm{1}}_{X_i\in k\mbox{\tiny -NN}(X)}{\mathbbm{1}}_{X_i\neq X} \\
& \leq \E \frac 1 k \sum \chi_U(X){\mathbbm{1}}_{X\in k\mbox{\tiny -NN}^{(X_1,X_2,\ldots,X_{i-1},X,X_{i+1},\ldots,X_n)}(X_i)}
{\mathbbm{1}}_{X_i\neq X} \\
 &\leq \E \chi_U(X)\frac 1 k (k C(d))\\
 & = C(d)\mu(U).
\end{align*}

Getting back to the term $(III)$, in the case $\Omega=\R^d$ we obtain:

\begin{align*}
(III)&\leq \E \left\vert\psi_{n,k}(X)-\eta_{n,k}(X)\right\vert 
\mathbbm{1}_{X\in K} + \E \left\vert\psi_{n,k}(X)-\eta_{n,k}(X)\right\vert 
\mathbbm{1}_{X\in U}  \\
& \leq C(d)\mu(U) + \mu(U) \\
&<(C(d)+1)\ve.
\end{align*}

Since $\ve>0$ is as small as desired, we conclude that $\eta_{n,k}(X)\to \eta(X)$ in probability, and so the $k$-NN classifier in $\ell^2(d)$ is universally (weakly) consistent in $\ell^2(d)$.

\section{Nagata dimension of a metric space}

Recall that a family $\gamma$ of subsets of a set $\Omega$ has multiplicity $\leq\delta$ if the intersection of more than $\delta$ different elements of $\gamma$ is always empty. In other words,
\[\forall x\in\Omega,~~\sum_{V\in\gamma}\chi_V(x)\leq\delta.\]

\begin{dfntn}
Let $\delta\in\N$, $s\in (0,+\infty]$. We say that a metric space $(\Omega,d)$ has {\em Nagata dimension} $\leq\delta$ {\em on the scale $s>0$}, if every finite family $\gamma$ of closed balls of radii $< s$ admits a subfamily $\gamma^\prime$ of multiplicity $\leq\delta+1$ which covers the centres of all the balls in $\gamma$. The smallest $\delta$ with this property, if it exists, and $+\infty$ otherwise, is called the Nagata dimension of $\Omega$ on the scale $s>0$.
A space $\Omega$ has Nagata dimension $\delta$ if it has Nagata dimension $\delta$ on a suitably small scale $s\in (0,+\infty]$. Notation: $\dim^s_{Nag}(\Omega)=\delta$, or simply $\dim_{Nag}(\Omega)=\delta$.
\end{dfntn}

Sometimes the following reformulation is more convenient.

\begin{prpstn} A metric space $(\Omega,d)$ has Nagata dimension $\leq\delta$ on the scale $s\in (0,+\infty]$ if and only if it satisfies the following property. For every $x\in\Omega$, $r<s$, and a sequence $x_1,\ldots,x_{\delta+2}\in \bar B_r(x)$, there are $i,j$, $i\neq j$, such that $d(x_i,x_j)\leq \max\{d(x,x_i),d(x,x_j)\}$.
\label{p:nagata2}
\end{prpstn}

\begin{proof} {\em Necessity} ($\Rightarrow$): from the family of closed balls $B_{d(x_i,x)}(x_i)$, $i=1,2,\ldots,\delta+2$, all having radii $<s$, extract a family of $\delta+1$ balls covering the centres. One of those balls, say with centre at $x_i$, must contain some $x_j$ with $i\neq j$, which means $d(x_i,x_j)\leq d(x_i,x)\leq \max\{d(x,x_i),d(x,x_j)\}$.

{\em Sufficiency} ($\Leftarrow$): let $\gamma$ be a finite family of closed balls of radii $< s$. Suppose it has multiplicity $>\delta+1$. Then there exist a point $x\in\Omega$ and $\delta+2$ balls in $\gamma$ with centres that we denote $x_1,\ldots,x_{\delta+2}$, all containing $x$. Denote $r=\max_i\{d(x,x_i)\}$. Then $r<s$, and by the hypothesis, there are $i,j$, $i\neq j$, with $d(x_i,x_j)\leq \max\{d(x,x_i),d(x,x_j)\}$. Without loss in generality, assume $d(x_i,x_j)\leq d(x_i,x)$, that is, $x_j$ belongs to the ball with centre in $x_i$. Now the ball centred at $x_j$ can be removed from the family $\gamma$, with the remaining family still covering all the centres and having the cardinality $\abs\gamma-1$. After finitely many steps, we arrive at a subfamily of multiplicity $\leq \delta+1$ covering all the centres.
\end{proof}

\begin{xmpl} 
The property of a metric space $\Omega$ having Nagata dimension zero on the scale $+\infty$ is equivalent to $\Omega$ being a non-archimedian metric space, that is, a metric space satisfying the strong triangle inequality, $d(x,z)\leq \max\{d(x,y),d(y,z)\}$.

Indeed, $\dim_{Nag}^{+\infty}(\Omega)=0$ means exactly that for any sequence of $\delta+2=2$ points, $x_1,x_2$, contained in a closed ball $B_r(x)$, we have $d(x_1,x_2)\leq \max\{d(x,x_1),d(x,x_2)\}$.
\label{ex:nonarchimedian}
\end{xmpl}

\begin{xmpl} 
It follows from Proposition \ref{p:nagata2} that $\dim_{Nag}(\R)=1$. Let $x_1,x_2,x_3$ be three points contained in a closed ball, that is, an interval $[x-r,x+r]$. Without loss in generality, assume $x_1<x_2<x_3$. If $x_2\leq x$, then $\abs{x_1-x_2}\leq \abs{x_1-x}$, and if $x_2\geq x$, then $\abs{x_3-x_2}\leq \abs{x_3-x}$.
\label{ex:R}
\end{xmpl}

The following example suggests that the Nagata dimension is relevant for the study of the $k$-NN classifier, as it captures in an abstract context the geometry behind Stone's lemma.

\begin{xmpl}
The Nagata dimension of the Euclidean space $\ell^2(d)$ is finite, and it is bounded by $C(d)-1$, where $C(d)$ is the value of the constant in Stone's geometric lemma (Lemma \ref{l:stonek=k}).

Indeed, let $x_1,\ldots,x_{C(d)+1}$ be points belonging to a ball with centre $x$. Using the argument in the proof of Stone's geometric lemma with $k=1$, mark $\leq C(d)$ points $x_i$ belonging to the $\leq C(d)$ cones with apex at $x$. At least one point, say $x_j$, has not been marked; it belongs to some cone, which therefore already contains a marked point, say $x_i$, different from $x_j$, and $\norm{x_i-x_j}\leq\norm{x_j-x}$.
\label{ex:cones}
\end{xmpl}

\begin{xmpl} 
A similar argument shows that every finite-dimensional normed space has finite Nagata dimension.
\end{xmpl}

\begin{rmrk}
In $\R^2=\C$ the family of closed balls of radius one centred at the vectors $\exp(2\pi ki/5)$, $k=1,2,\ldots,5$, has multiplicity $5$ and admits no proper subfamily containing all the centres. Therefore, the Nagata dimension of $\ell^2(2)$ is at least $5$. Since the plane can be covered with $6$ cones having the central angle $\pi/3$, Example \ref{ex:cones} implies that $\dim_{Nag}(\ell^2(2))= 5$. 
\end{rmrk}

\begin{rmrk}
The problem of calculating the Nagata dimension of the Euclidean space $\ell^2(d)$ is mentioned as ``possibly open'' by Nagata \cite{nagata_open}, p. 9 (where the value $\dim_{Nag}+2$ is called the ``crowding number''). Nagata also remarks that $\dim_{Nag}\R^1=1$ and $\dim_{Nag}\ell^2(3)=5$ (without a proof).
\end{rmrk}

\begin{rmrk}
Notice that the property of the Euclidean space established in the proof of Stone's geometric lemma is strictly stronger than the finiteness of the Nagata dimension. There exists a finite $\delta$ (in general, higher than the Nagata dimension) such that, given a sequence $x_1,\ldots,x_{\delta+2}\in \bar B_r(x)$, $r<s$, there are $i,j$, $i\neq j$, such that $d(x_i,x_j)< \max\{d(x,x_i),d(x,x_j)\}$. The inequality here is strict, cf. Remark \ref{r:strictineq}. This is exactly the property that removes the problem of distance ties in the Euclidean space. However, adopting this as a definition in the general case would be too restrictive, removing from consideration a large class of metric spaces in which the $k$-NN classifier is still universally consistent, such as all non-archimedean metric spaces.
\end{rmrk}

\begin{xmpl} Let $e_n$ denote the $n$-th standard basic vector in the separable Hilbert space $\ell^2$, that is, a sequence whose $n$-th coordinate is $1$ and the rest are zeros.
The convergent sequence $(1/n)e_n$, $n\geq 0$, together with the limit $0$, viewed as a metric subspace of $\ell^2$, has infinite Nagata dimension on every scale $s>0$. This is witnessed by the family of closed balls $B_{1/n}((a/n)e_n)$, having zero as the common point, and having the property that every centre belongs to exactly one ball of the family. Realizing $\R$ as a continuous curve in $\ell^2$ without self-intersections passing through all elements of the sequence as well as the limit leads to an equivalent metric on $\R$ having infinite Nagata dimension on each scale.
\end{xmpl}

\begin{rmrk}
The Nagata--Ostrand theorem \cite{nagata,ostrand} states that the Lebesgue covering dimension of a metrizable topological space is the smallest Nagata dimension of a compatible metric on the space (and in fact this is true on every scale $s>0$, \cite{nagata1}). This is the historical origin of the concept of the metric dimension.

There appears to be no single comprehensive reference to the concept of Nagata dimension. Various results are scattered in the journal papers \cite{AQ,nagata1,nagata,nagata_open,ostrand,preiss83}, see also the book \cite{nagata_dim_theory}, pages 151--154.
\end{rmrk}

Metric spaces of finite Nagata dimension admit an almost literal version of Stone's geometric lemma in case where the sample has no distance ties, that is, the values of the distances $d(x_i,x_j)$, $i\neq j$, are all pairwise distinct.

\begin{lmm}[Stone's geometric lemma, finite Nagata dimension, no ties]
Let $\Omega$ a metric space of Nagata dimension $\delta<\infty$ on a scale $s>0$.
Let 
\[\sigma=(x_1,x_2,\ldots,x_n),~i=1,2,\ldots,n,\]
be a finite sample in $\Omega$, and let $x\in \Omega$ be any. Suppose there are no distance ties inside the sample
\begin{equation}
x,x_1,x_2,\ldots, x_{i-1},x_i,x_{i+1},\dots,x_n,
\label{eq:samplexxi}
\end{equation}
and $k$ is such that, inside the above sample, $r_{k\mbox{\tiny -NN}}(x_i)<s$ for all $i$. The number of $i$ having the property that $x\neq x_i$ and $x$ is among the $k$ nearest neighbours of $x_i$ inside the sample above is limited by $(k+1)(\delta+1)$. 
\label{l:stone_sem_empates}
\end{lmm}

\begin{proof}
Suppose that $i=1,2,\ldots,m$ is such that $x_i\neq x$ and $x_i$ has $x$ among the $k$ nearest neighbours inside the sample as in Eq. (\ref{eq:samplexxi}).
The family $\gamma$ of closed balls $B_{r_{k\mbox{\tiny -NN}}(x_i)}(x_i)$, $i\leq m$, admits a subfamily $\gamma^\prime$ of multiplicity $\leq\delta+1$ covering all the points $x_i$, $i\leq m$. Since there are no distance ties, every ball belonging to $\gamma$ contains $\leq k+1$ points. It follows that $\sharp\gamma^\prime\geq m/(k+1)$. All the balls in $\gamma^\prime$ contain $x$, and we conclude: $\sharp\gamma^\prime\leq\delta+1$. The result follows.
\end{proof}

Now the same argument as in the original proof of Stone (Subs. \ref{ss:stone}) shows that the $k$-NN classifier is consistent under each distritution $\mu$ on $\Omega\times\{0,1\}$ with the property that the distance ties occur with zero probability. Since we are going to give a proof of a more general result, we will not repeat the argument here, only mention that due to the Cover--Hart lemma, if $n$ is sufficiently large, then with arbitrarily high probability, the $k$ nearest neighbours of a random point inside a random sample will all lie at a distance $<s$.

\section{Distance ties}

In this section we will construct a series of examples to illustrate the difficulties arising in the presence of distance ties in general metric spaces that are absent in the Euclidean case.
The fundamental difference between the two situations is the inequality in the equivalent definition of the Nagata dimension (proposition \ref{p:nagata2}) that is, unlike in the Euclidean space, no longer strict. 

As we have already noted (remark \ref{r:unimportant}), the conclusion of  Stone's geometric lemma \ref{l:stonek=k} remains valid even if we allow the adversary to break the distance ties and pick the $k$ nearest neighbours. Our first example shows that it is no longer the case in a metric space of finite Nagata dimension.

\begin{xmpl}
Consider a finite set $\sigma=\{x_1,x_2,\ldots,x_n\}$ with $n\geq k$ points, and assume that in the metric space $\sigma\cup\{x\}$ all $n+1$ points are pairwise at a distance one from each other. The Nagata dimension of the metric space $\sigma\cup\{x\}$ is equal to $\delta=0$. Indeed, if a family $\gamma$ of closed balls contains any ball of radius $\geq 1$, it already covers $\sigma$ on its own. Otherwise, we choose one ball of radius $<1$ (that is, a singleton) for each centre. The multiplicity of the selected subfamily is $1$ in each case.

Now let us discuss the distance ties.
For any element $x_i$ of $\sigma$, the remaining $n$ points of $\sigma\cup\{x\}$ are tied between themselves as the possible $k$ nearest neighbours. The adversary may decide to always select $x$ among them, thus invalidating the conclusion of Stone's geometric lemma.

However, the problem is easily resolved if we break distance ties using a uniform distribution on the nearest neighbour candidates. In this case, the expected number of indices $i$ such that $x$ is chosen as one of the $k$ nearest neighbours of $x_i$ within the sample $\{x_1,x_2,\ldots,x,\ldots,x_{k}\}$ is obviously $k$.
\end{xmpl}

\begin{rmrk}
It is worth observing that in the Euclidean case $\Omega=\R^d$ the size of a sample inheriting a $0$-$1$ distance will be limited from above by the dimension, $d$. 
\end{rmrk}

The next example shows that Stone's geometric lemma in finite dimensional metric spaces cannot be saved even with the uniform tie-breaking.

\begin{xmpl} There exists a countable metric space $\sigma=\{x_1,x_2,\ldots\}$ of Nagata dimension $0$, having the following property.
Given $N\in\N$, for a sufficiently large $n$ the expected number of points $x_i\neq x_1$ within the sample $\sigma_n=\{x_1,\ldots,x_n\}$ having $x_1$ as the nearest neighbour under the uniform tie-breaking is $\geq N$.

We will construct $\sigma$ recurrently. Let $\sigma_1=\{x_1\}$. Add $x_2$ at a distance $1$ from $x_1$, and set $\sigma_2=\{x_1,x_2\}$. If $\sigma_n$ has been already defined, add $x_{n+1}$ at a distance $2^n$ from all the existing points $x_i$, $i\leq n$, and set $\sigma_{n+1}=\sigma_n\cup\{x_{n+1}\}$. It is clear that the distance so defined is a metric. 

We will verify by induction in $n$ that $\dim_{Nag}(\sigma_n)=0$, on the scale $s=+\infty$. For $n=1$ this is trivially true. Assume the statement holds for $\sigma_n$, and let $\gamma$ be a family of closed balls in $\sigma_{n+1}$. If one of those balls contains all the points, there is nothing to prove. Assume not, that is, all the balls elements of $\gamma$ have radii smaller than $2^{n+1}$. Choose a subfamily of multiplicity $1$ consisting of balls centred in elements of $\sigma_n$ and covering them all, and add one ball centred in $x_{n+1}$ (which is a singleton). Now it follows that $\dim_{Nag}\sigma=0$ as well.

Finally, let us show that if $n$ is sufficiently large, then the expected number of indices $i$ such that $x=x_1$ is the nearest neighbour of $x_i$ under a uniform tie-breaking is as large as desired. With this purpose, for each $i\geq 2$ we will calculate the expectation of the event $x_1\in NN(x_i)$, where $NN(x_i)$ denotes the set of nearest neighbours of $x_i$ in the rest of the finite sample $\sigma_n$.

For $x_2$, the unique nearest neighbour within $\sigma_n$ is $x_1$, therefore $\E[x_1\in NN(x_2)]=1$. For $x_3$, there are two points in $\sigma_n$ at a distance $2$ from $x_3$, which can be chosen each with probability $1/2$, namely $x_1$ and $x_2$, therefore $\E[x_1\in NN(x_3)]=1/2$. For arbitrary $i$, in a similar way, $\E[x_1\in NN(x_i)]=1/i$. We conclude:
\[\E[\sharp\{i=1,\ldots,n\colon x_1\in NN(x_i)\}=\sum_{i=1}^n\frac 1i,\]
and the sum of the harmonic series converges to $+\infty$ as $n\to\infty$.
\label{ex:exemplonagata}
\end{xmpl}

Can it be that the distance ties are in some sense extremely rare? Even this expectation is unfounded. 

\begin{xmpl}
Given a value $\delta>0$ (risk) and a sequence $n^\prime_k\uparrow+\infty$, there exist a compact metric space of Nagata dimension zero (a Cantor space  with a suitable compatible metric) equipped with a non-atomic probability measure, and a sequence $n_k\geq n^{\prime}_k$, $k/n_k\to 0$, with the following property. With confidence $>1-\delta$, for every $k$, a random element $X$ has $\geq n_k$ distance ties among its $k$ nearest neighbours within a random $n_{k+1}$-sample $\sigma$. 

The space $\Omega$, just like in the Preiss example (Sect. \ref{s:preiss}), is the direct product $\prod_{k=1}^{\infty} [N_k]$ of finite discrete spaces, whose cardinalities $N_k\geq 2$ will be chosen recursively, and $[N_k]=\{1,2,\ldots,N_k\}$. The metric is given by the rule
\[d(\sigma,\tau)=\begin{cases} 0,&\mbox{ if }\sigma=\tau,\\
2^{-\min\{i\colon \sigma_i\neq\tau_i\}},&\mbox{ otherwise.}
\end{cases}\]
This metric induces the product topology and is non-archimedian, so the Nagata dimension of $\Omega$ is zero (example \ref{ex:nonarchimedian}). The measure $\mu$ is the product of uniform measures $\mu_{N_k}$ on the spaces $[N_k]$. This measure is non-atomic, and in particular, $\mu$-almost all distance ties occur at a strictly positive distance from a random element $X$.

Choose a sequence $(\delta_i)$ with $\delta_i>0$ and $2\sum_i\delta_i=\delta$. Choose $N_1$ so large that, with probability $>1-\delta_1$, $n_1=n^\prime_1$ independent random elements following a uniform distribution on the space $[N_1]$ are pairwise distinct. Now let $n_2\geq n^\prime_2$ be so large that with probability $>1-\delta_1$, if $n_2$ independent random elements follow a uniform distribution on $[N_1]$, then each element of $[N_1]$ appears among them at least $n_1$ times. 

Suppose that $n_1,N_1,n_2,N_2,\ldots,n_k$ have been chosen. Let $N_{k}$ be so large that, with probability $>1-\delta_{k}$, $n_k$ i.i.d. random elements uniformly distributed in $[N_k]$ are pairwise distinct. Choose $n_{k+1}\geq n^\prime_{k+1}$ so large that, with probability $>1-\delta_{k}$, if $n_{k+1}$ i.i.d. random elements are uniformly distributed within $\prod_{i=1}^k [N_i]$, then each element of $\prod_{i=1}^k [N_i]$ will appear among them at least $n_k$ times.

Let $k$ be any positive natural number. Choose $n_{k+1}+1$ i.i.d. random elements $X,X_1,\ldots,X_{n_{k+1}}$ of $\Omega$, following the distribution $\mu$. With probability $>1-2\delta_k$, the following occurs: there are $n_k$ elements in the sample $X_1,X_2,\ldots,X_{n_k}$ which have the same $i$-th coordenates as $X$, $i=1,2,3,\ldots,k$, yet the $(k+1)$-coordenates of $X,X_1,\ldots,X_{n_k}$ are all pairwise distinct. In this way, the distances between $X$ and all those $n_k$ elements are equal to $2^{-k-1}$. We have $n_k$ distance ties between $k$ nearest neighbours of $X$ (which are all at the same distance as the nearest neighbour of $X$), and $n_k\geq n_k^\prime$, as desired.
\end{xmpl}

Now, it would be tempting to try and reduce the general case to the case of zero probability of ties, as follows. 
Recall that the $\ell^1$-type direct sum of two metric spaces, $(X,d_X)$ and $(Y,d_Y)$, is the direct product $X\times Y$ equipped with the coordinatewise sum of the two metrics:
\[d(x,y)=d_X(x_1,x_2)+d_Y(y_1,y_2).\]
Notation: $X\oplus_1 Y$.

Let $\Omega$ be a domain, that is, a metric space equipped with a probability measure $\mu$ and a regression function, $\eta$. Form the $\ell^1$-type direct sum $\Omega\oplus_1 [0,\ve]$, and equip it with the product measure $\mu\otimes\lambda$ (where $\lambda$ is the normalized Lebesgue measure on the interval) and the regression function $\eta\circ\pi_1$, where $\pi_1$ is the production on the first coordinate. It is easy to see that the probability of distance ties in the space $\Omega\oplus_1 [0,\ve]$ is zero, and every uniform distance tie breaking within a given finite sample will occur for a suitably small $\ve>0$. In this way, one could derive the consistency of the classifier by conitioning. However, we will now give an example of two metric spaces of Nagata dimension $0$ and $1$ respectively, whose $\ell^1$-type sum has infinite Nagata dimension. This is again very different from what happens in the Euclidean case. 

\begin{xmpl}
Fix $\alpha>0$.
Let $\Omega=\{x_n\colon n\in\N\}$, equipped with the following distance: 
\[d(x_i,x_j) = \begin{cases} 0,&\mbox{ if }i=j,\\
\sum_{k=1}^j \alpha^k,&\mbox{ if }i<j.
\end{cases}\] 
For $i<j<k$,
\[d(x_i,x_k) = \sum_{m=1}^k \alpha^m = d(x_j,x_k) > d(x_i,x_j) = \sum_{m=1}^j \alpha^m,\]
from where it follows that $d$ is an ultrametric. Thus, $\Omega$ is a metric space of Nagata dimension $0$.

The interval $\I=[0,1]$ has Nagata dimension $1$ (Ex. \ref{ex:R}). Now let us consider the $\ell^1$-type sum $\Omega \oplus_1 \I$. Let $0<\beta<\alpha<1$, and $\beta<1/2$.
Consider the infinite sequence
\[z_i = (x_i,\beta^{i})\in \Omega \oplus_1 \I\]
and the point 
\[z = (x_0,0).\]
Whenever $i<j$,
\begin{align*}
d(z_i,z_j) & = d(x_i,x_j) + \beta^i-\beta^j \\
& \geq d(x_i,x_0) + \alpha^j+ \beta^i-\beta^j \\
& > d(x_i,x_0) + \beta^i \\
& = d(z_i,z),
\end{align*}
and also
\begin{align*}
d(z_i,z_j) & = d(x_i,x_j) + \beta^i-\beta^j \\
& > d(x_j,x_0) + \beta^j \\
& = d(z_j,z).
\end{align*}
Together, the properties imply: for all $i\neq j$,
\[d(z_i,z_j) > \max\{d(z_i,z_0),d(z_j,z_0) \}.\]
Thus, the Nagata dimension of the $\ell^1$-type sum $\Omega \oplus_1 \I$ is infinite.
\end{xmpl}

The above examples show that beyond the Euclidean setting, we have to put up with the possibility that some points in a sample will appear disproportionately often among $k$ nearest neighbours of other points. In data science, such points are known as ``hubs'' and the above (empirical) observation, as the ``hubness phenomenon'',  see e.g., \cite{RNI} and further references therein. Stone's geometric lemma has to be generalized to allow for the possibility of a few of those ``hubs'', whose number will be nevertheless limited. The lemma has to be reshaped in the spirit of the Hardy--Littlewood inequality in geometric measure theory.

To begin with, following Preiss \cite{preiss83}, we will extend further our metric space dimension theory setting.

\section{Sigma-finite dimensional metric spaces}

\begin{dfntn} Say that a metric subspace $X$ of a metric space $\Omega$ has {\em Nagata dimension $\leq\delta\in\N$ on the scale $s>0$ inside of} $\Omega$ if every finite family of closed balls in $\Omega$ with centres in $X$ admits a subfamily of multiplicity $\leq\delta+1$ in $\Omega$ which covers all the centres of the original balls. The subspace $X$ has a finite Nagata dimension in $\Omega$ if $X$ has finite dimension in $\Omega$ on some scale $s>0$. Notation: $\dim^s_{Nag}(X,\Omega)$ or sometimes simply $\dim_{Nag}(X,\Omega)$.
\label{d:dimnagata}
\end{dfntn}

Following Preiss, let us call a family of balls {\em disconnected} if the centre of each ball does not belong to any other ball. Here is a mere reformulation of the above definition.

\begin{prpstn}  
For a subspace $X$ of a metric space $\Omega$, one has
\[\dim^s_{Nag}(X,\Omega) \leq \delta\]
if and only if every disconnected family of closed balls in $\Omega$ of radii $<s$ with centres in $X$ has multiplicity $\leq \delta+1$.
\label{ex:famdesconexa}
\end{prpstn}

\begin{proof}
{\em Necessity.} Let $\gamma$ be a disconnected finite family of closed balls in $\Omega$ with centres in $X$. Since by assumption $\dim^s_{Nag}(X,\Omega) \leq \delta$, the family $\gamma$ admits a subfamily of multiplicity $\leq\delta+1$ covering all the original centres. But only subfamily that contains all the centres is $\gamma$ itself.

{\em Sufficiency.} 
Let $\gamma$ be a finite family of closed balls in $\Omega$ with centres in $X$.
Denote $C$ the set of centres of those balls. 
Among all the disconnected subfamilies of $\gamma$ (which exist, e.g., each family containing just one ball is such) there is one, $\gamma^\prime$, with the maximal cardinality of the set $C\cap \cup\gamma^\prime$. We claim that $C\subseteq\gamma^\prime$, which will finish the argument. Indeed, if it is not the case, there is a ball, $B\in\gamma$, whose centre, $c\in C$, does not belong to $\cup \gamma^\prime$. Remove from $\gamma^\prime$ all the balls with centres in $B\cap C$ and add $B$ instead. The new family, $\gamma^{\prime\prime}$, is disconnected and contains $\left(C\cap\cup\gamma^\prime\right)\cup\{c\}$, which contradicts the maximality of $\gamma^\prime$. 
\end{proof}

In the definition \ref{d:dimnagata}, as well as in the proposition \ref{ex:famdesconexa}, closed balls can be replaced with open ones. In fact, the statements remain valid if some balls in the families are allowed to be closed, other, open. We have the following.

\begin{prpstn} 
For a subspace $X$ of a metric space $\Omega$, the following are equivalent.
\begin{enumerate}
\item\label{equiv:1} $\dim^s_{Nag}(X,\Omega) \leq\delta$,
\item\label{equiv:3} every finite family of balls (some open, others closed) in $\Omega$ with centres in $X$ and radii $<s$ admits a subfamily of multiplicity $\leq\delta+1$ in $\Omega$ which covers all the centres of the original balls,
\item\label{equiv:2} every finite family of open balls in $\Omega$ having radii $<s$ with centres in $X$ admits a subfamily of multiplicity $\leq\delta+1$ in $\Omega$ which covers all the centres of the original balls,
\item\label{equiv:4} every disconnected family of open balls in $\Omega$ of radii $<s$ with centres in $X$ has multiplicity $\leq\delta+1$,
\item\label{equiv:5} every disconnected family of balls (some open, others closed) in $\Omega$ of radii $<s$ with centres in $X$ has multiplicity $\leq\delta+1$.
\end{enumerate}
\end{prpstn}

\begin{proof}
$(\ref{equiv:1})\Rightarrow(\ref{equiv:3})$: Let $\gamma$ be a finite family of  balls in $\Omega$ with centres in $X$, of radii $<s$, where some of the balls may be open and others, closed. For every element $B\in\gamma$ and each $k\geq 2$, form a closed ball $B_k$ as follows: if $B$ is closed, then $B_k=B$, and if $B$ is open, then define $B_k$ as having the same centre and radius $r(1-1/m)$, where $r$ is the radius of $B$. Thus, we always have $B=\cup_{k=2}^{\infty}B_k$. Select recursively a chain of subfamilies 
\[\gamma\supseteq \gamma_1\supseteq\gamma_2\supseteq \ldots \supseteq \gamma_k\supseteq \ldots\]
with the properties that for each $k$, the family of closed balls $B_k$, $B\in\gamma_k$ has multiplicity $\leq\delta+1$ in $\Omega$ and covers all the centres of the balls in $\gamma$. Since $\gamma$ is finite, starting with some $k$, the subfamily $\gamma_k$ stabilizes, and now it is easy to see that the subfamily $\gamma_k$ itself has the desired multiplicity, and of course covers all the original centres.

$(\ref{equiv:3})\Rightarrow(\ref{equiv:2})$: Trivially true.

$(\ref{equiv:2})\Rightarrow(\ref{equiv:4})$:
Same argument as in the proof of necessity in proposition \ref{ex:famdesconexa}.

$(\ref{equiv:4})\Rightarrow(\ref{equiv:5})$: Let $\gamma$ be a disconnected family of balls in $\Omega$, some of which may be open and others, closed, having radii $<s$ and centred in $X$. For each $B\in\gamma$ and $\ve>0$, denote $B_{\ve}$ an open ball equal to $B$ if $B$ is open, and concentric with $B$ and of the radius $r+\ve$, where $r$ is the radius of $B$, if $B$ is closed.
For a sufficiently small $\ve>0$, the family $\{B_\ve\colon B\in\gamma\}$ is disconnected, and its radii are all strictly less than $s$, therefore this family has multiplicity $\leq\delta+1$ by assumption. The same follows for $\gamma$.

$(\ref{equiv:5})\Rightarrow(\ref{equiv:1})$: the condition $(\ref{equiv:5})$ is formally even stronger than an equivalent condition for $(\ref{equiv:1})$ established in Proposition \ref{ex:famdesconexa}.
\end{proof}

\begin{prpstn}
Let $X$ be a subspace of a metric space $\Omega$, satisfying $\dim^s_{Nag}(X,\Omega)\leq \delta$. Then $\dim^s_{Nag}(\bar X,\Omega)\leq \delta$, where $\bar X$ is the closure of $X$ in $\Omega$.
\label{p:closurenagata}
\end{prpstn}

\begin{proof}
Let $\gamma$ be a finite disconnected family of open balls in $\Omega$ of radii $<s$, centred in $\bar X$. Let $y\in\Omega$, and let $\gamma^\prime$ consist of all balls in $\gamma$ containing $y$. Choose $\ve>0$ so small that the open $\ve$-ball around $y$ is contained in every element of $\gamma^\prime$. For every open ball $B\in\gamma^\prime$, denote $y_B$ the centre and $r_B$ the radius. We can also assume that $\ve<r_B$ for each $B\in\gamma^\prime$. Denote $B^\prime$ an open ball of radius $r_B-\ve>0$, centred at a point $x_B\in X$ satisfying $d(x_B,y_B)<\ve$. Then $B^\prime\subseteq B$, so the family $\{B^\prime\colon B\in\gamma^\prime\}$ is disconnected, and has radii $<s$. Therefore, $y$ only belongs to $\leq\delta+1$ balls $B^\prime$, $B\in\gamma^\prime$, consequently the cardinality of $\gamma^\prime$ is bounded by $\delta+1$.
\end{proof}

\begin{prpstn}
If $X$ and $Y$ are two subspaces of a metric space $\Omega$, having finite Nagata dimension in $\Omega$ on the scales $s_1$ and $s_2$ respectively, then $X\cup Y$ has a finite Nagata dimension in $\Omega$, with $\dim_{Nag}(X\cup Y,\Omega)\leq \dim_{Nag}(X,\Omega)+\dim_{Nag}(Y,\Omega)$, on the scale $\min\{s_1,s_2\}$.
\label{p:uniaonagata}
\end{prpstn}

\begin{proof}
Given a finite family of balls $\gamma$ in $\Omega$ of radii $<\min\{s_1,s_2\}$ centred in elements of $X\cup Y$, represent it as $\gamma=\gamma_X\cup\gamma_Y$, where the balls in $\gamma_X$ are centered in $X$, and the balls in $\gamma_Y$ are centred in $Y$. The rest is obvious.
\end{proof}

\begin{dfntn}
A metric space $\Omega$ is said to be {\em sigma-finite dimensional in the sense of Nagata} if $\Omega=\cup_{i=1}^{\infty}X_n$, where every subspace $X_n$ has finite Nagata dimension in $\Omega$ on some scale $s_n>0$ (where the scales $s_n$ are possibly all different).
\end{dfntn}

\begin{rmrk}
Due to Proposition \ref{p:closurenagata}, in the above definition we can assume the subspaces $X_n$ to be closed in $\Omega$, in particular Borel subsets.
\end{rmrk}

\begin{rmrk}
A good reference for a great variety of metric dimensions, including the Nagata dimension, and their applications to measure differentiation theorems, is the article \cite{AQ}.
\end{rmrk}

Now we will develop a version of Stone's geometric lemma for general metric spaces of finite Nagata dimension.

\section{From Stone to Hardy--Littlewood} 

\begin{lmm}
Let $\sigma=\{x_1,x_2,\ldots,x_n\}$ be a finite sample in a metric space $\Omega$, and let $X$ be a subspace of finite Nagata dimension $\delta$ in $\Omega$ on a scale $s>0$. Let $\alpha\in (0,1]$ be any. Let $\sigma^\prime\sqsubseteq\sigma $ be a sub-sample with $m$ points. Assign to every $x_i\in\sigma$ a ball, $B_i$ (which could be open or closed), centred at $x_i$, of radius $<s$. Then 
\[\sharp\{i=1,2,\ldots,n\colon x_i\in X,~~ \sharp(B_i\cap \sigma^\prime)\geq\alpha\sharp (B_i\cap\sigma)\}\leq \alpha^{-1}(\delta+1)m.
\]
\label{l:alphadelta}
\end{lmm}

\begin{proof}
Denote $I$ the set of all $i=1,2,\ldots,n$ such that 
\[\sharp(B_i\cap \sigma^\prime)\geq\alpha\sharp (B_i\cap\sigma).\]
According to the assumption $\dim_{Nag}^s(\Omega)=\delta$, there exists $J\subseteq I$ such that the subfamily $\{B_i\colon i\in J\}$ has multiplicity $\leq\delta+1$ and covers all the centres $x_i$, $i\in I$. In particular, each point of $\sigma^\prime$ belongs, at most, to $\delta+1$ balls $B_i$, $i\in I$. Consequently,
\begin{align*}
\alpha\sharp I &= \alpha\sharp\{x_i\colon i\in I\} \\
&\leq \sum_{i\in I}\alpha\sharp(B_i\cap\sigma) \\
& \leq\sum_{i\in I}\sharp(B_i\cap\sigma^\prime) 
\\
 & \leq \sharp\sigma^\prime (\delta+1) \\
&= m(\delta+1),
\end{align*}
whence the conclusion follows.
\end{proof}

\begin{rmrk}
In applications of the lemma, $B_i=B_{\e_{\mbox{\tiny{k-NN}}}}(x)$ (sometimes the ball will need to be open, sometimes closed, depending on the presence of distance ties).
\end{rmrk}

\begin{lmm}
Let $\alpha,\alpha_1,\alpha_2\geq 0$, $t_1,t_2\in [0,1]$, $t_2\leq 1-t_1$. Assume that $\alpha_1\leq\alpha$ and
\[t_1\alpha_1+(1-t_1)\alpha_2\leq\alpha.\]
Then
\[\frac{t_1\alpha_1+t_2\alpha_2}{t_1+t_2}\leq\alpha.\]
\label{l:fracoes}
\end{lmm}

\begin{proof} If $\alpha_2\leq\alpha$, the conclusion is immediate. Otherwise, $\alpha_2>\alpha$, and it follows that
\[t_1\alpha_1+t_2\alpha_2\leq \alpha - (1-t_1-t_2)\alpha_2\leq (t_1+t_2)\alpha. \]
\end{proof}

\begin{lmm}
Let $x,x_1,x_2,\ldots,x_n$ be a finite sample (possibly with repetitions), and let $\sigma^\prime\sqsubset \sigma$ be a subsample. Let $\alpha\geq 0$, and let $B$ be a closed ball around $x$ of radius $r_{k\mbox{\tiny -NN}}(x)$ which contains $K$ elements of the sample,
\[\sharp\{i=1,2,\ldots,n\colon x_i\in B\}=K.\]
Suppose that the fraction of points of $\sigma^\prime$ found in $B$ is no more than $\alpha$, 
\[\sharp\{i\colon x_i\in\sigma^\prime,~~x_i\in B\}\leq \alpha K,\]
and that the same holds for the corresponding open ball, $B^\circ$,
\[\sharp\{i\colon x_i\in\sigma^\prime,~~x_i\in B^\circ\}\leq \alpha \sharp\{i\colon x_i\in B^\circ\}.\]
Under the uniform tie-breaking of the $k$ nearest neighbours, the expected fraction of the points of $\sigma^\prime$ found among the $k$ nearest neighbours of $x$ is less than or equal to $\alpha$.
\label{l:alphak}
\end{lmm}

\begin{proof}
We apply lemma \ref{l:fracoes} with $\alpha_1$ and $\alpha_2$ being the fractions of the points of $\sigma^\prime$ found in the closed ball $B$ and on the sphere $S=B\setminus B^\circ$ respectively, $t_1=\sharp B^\circ /\sharp B$, and $t_2$ being the fraction of the points of the sphere $S=B\setminus B^\circ$ to be chosen uniformly and randomly as the $k$ nearest neighbours of $x$ that are still missing in the open ball $B^\circ$. Now it is enough to observe that the expected fraction of the points of $\sigma^\prime$ among the $k$ nearest neighbours that belong to the sphere is also equal to $\alpha_2$, because they are being chosen randomly, following a uniform distribution.
\end{proof}

Now we can give a promised alternative proof of the principal result along the same lines as Stone's original proof in the finite-dimensional Euclidean case.

\begin{thrm}
The $k$ nearest neighbour classifier under the uniform distance tie-breaking is universally consistent in every metric space having sigma-finite Nagata dimension, when $n,k\to\infty$ and $k/n\to 0$.
\label{t:nagatasigmafinitaconsistente}
\end{thrm}

\begin{proof}
Represent $\Omega=\cup_{l=1}^\infty Y_n$, where $Y_n$ have finite Nagata dimension in $\Omega$. According to proposition \ref{p:uniaonagata}, we can assume that $Y_n$ form an increasing chain, and proposition \ref{p:closurenagata} allows to assume that $Y_n$ are Borel sets. Let a Borel probability measure $\mu$ and a measurable regression function $\eta$ be any on $\Omega$.
Given $\e>0$, there exists $l$ such that $\mu(Y_l)\geq 1-\e/2$. According to Luzin's theorem (Th. \ref{th:luzin}), there is a closed precompact subset $K\subseteq Y_l$, such that $\eta\vert_K$ is continuous and $\mu(K)\geq 1-\e$. Applying the Tietze--Urysohn extension theorem (\cite{E}, 2.1.8), we extend the function $\eta\vert_K$ to a continuous function $g$ over $\Omega$. 

In the spirit of the proof of Stone's theorem \ref{t:stone}, it is enough to limit the term
\begin{align*}
(B)  &=  \E \frac 1 k \sum_{i=1}^n \left\vert\eta(X_i)-g(X_i)\right\vert {\mathbbm 1}_{X_i\in k\mbox{\tiny -NN}(X)}{\mathbbm{1}}_{X\in K}{\mathbbm{1}}_{X_i\notin K}
 \\
&= \E  \E_{j\sim \mu_{\sharp}} \frac 1 k 
\sum_{i\in \{0,1,\ldots,n\}\setminus\{j\}} 
\left\vert\eta(X_i)-g(X_i)\right\vert 
{\mathbbm{1}}_{X_{i}\in k\mbox{\tiny -NN}(X_j)}{\mathbbm{1}}_{X_{j}\in K}{\mathbbm{1}}_{X_{i}\notin K},
\end{align*}
where $\mu_{\sharp}$ is the uniform measure on the set $\{0,1,2,\ldots,n\}$, and we denote $X_0=X$. 

Let $s>0$ denote the scale on which $Y_l$ has finite Nagata dimension, $\delta\in\N$.
Since by the Cover--Hart lemma \ref{l:cover-hart} 
\[\E_{j\sim \mu_{\sharp}}r^{\{X_0,\ldots,X_n\}\setminus \{X_j\}}_{k\mbox{\tiny -NN}}(X_j)
 \to 0,
\]
it suffices to estimate the term
\[(B^\prime) = \E  \E_{j\sim \mu_{\sharp}} \frac 1 k 
\sum_{i\in \{0,1,\ldots,n\}\setminus\{j\}} 
\left\vert\eta(X_i)-g(X_i)\right\vert 
{\mathbbm{1}}_{X_{i}\in k\mbox{\tiny -NN}(X_j)}{\mathbbm{1}}_{X_{j}\in K}{\mathbbm{1}}_{X_{i}\notin K}{\mathbbm{1}}_{r^{\{X_0,\ldots,X_n\}\setminus \{X_j\}}_{k\mbox{\tiny -NN}}(X_j)<s}.
\]

We will treat the above as a sum of two expectations, $(B_1)$ and $(B_2)$, according to whether the $k$ nearest neighbours of $X_j$ inside the sample $\{X_0,X_1,\ldots,X_{j-1},X,X_{j+1},\ldots X_n\}$ contain more or less than $\sqrt\e k$ elements belonging to $U=\Omega\setminus K$. 

Applying lemma \ref{l:alphadelta} to the closed balls of radius $k\mbox{-NN}(X_j)$ as well as the corresponding open balls, together with lemma \ref{l:alphak}, we get in the first case
\begin{align*}
(B_1) &= \E  \E_{j\sim \mu_{\sharp}} \frac 1 k \sum \left\vert\eta(X_i)-g(X_i)\right\vert 
{\mathbbm{1}}_{X_{i}\in k\mbox{\tiny -NN}(X_j)}{\mathbbm{1}}_{X_{j}\in K}{\mathbbm{1}}_{X_{i}\notin K,~ i\in \{0,1,\ldots,n\}\setminus\{j\}} {\mathbbm{1}}_{X_{i}\notin K}{\mathbbm{1}}_{r^{\{X_0,\ldots,X_n\}\setminus \{X_j\}}_{k\mbox{\tiny -NN}}(X_j)<s}\times
\\ 
& \times {\mathbbm{1}}_{\sharp\{i\colon X_i\in k\mbox{\tiny -NN}(X_j),~~X_i\notin K\}\geq k\sqrt \e}
 \\
&\leq 
\E  \frac 1 k k \e^{-1/2} (\delta+1) \frac 1n\sharp\{i=0,1,\ldots,n\colon X_i\notin K\}\\
&\leq  \e^{-1/2} (\delta+1) \e\\
& = \sqrt{\e}(\delta+1),
\end{align*}
where we have used the fact that the sum in the first line does not exceed $k$.
In the second case,
\begin{align*}
(B_2) &= \E  \E_{j\sim \mu_{\sharp}} \frac 1 k \sum \left\vert\eta(X_i)-g(X_i)\right\vert 
{\mathbbm{1}}_{X_{i}\in k\mbox{\tiny -NN}(X_j)}{\mathbbm{1}}_{X_{j}\in K}{\mathbbm{1}}_{X_{i}\notin K,~ i\in \{0,1,\ldots,n\}\setminus\{j\}}{\mathbbm{1}}_{r^{\{X_0,\ldots,X_n\}\setminus \{X_j\}}_{k\mbox{\tiny -NN}}(X_j)<s}\times
\\ 
& \times {\mathbbm{1}}_{\sharp\{i\colon X_i\in k\mbox{\tiny -NN}(X_j),~~X_i\notin K\}< k\sqrt \e}
\\
&\leq \frac 1kk\sqrt \e \\
& = \sqrt \e.
\end{align*}
\end{proof}

\renewcommand\thesection{\Alph{section}}
\setcounter{section}{1}
\section*{Appendix: Luzin's theorem}

The classical Luzin theorem admits numerous variations, of which we need the following one.

\begin{thrm}[Luzin's theorem] Let $X$ be a separable metric space (not necessarily complete), $\mu$ a Borel probability measure on $X$, and $f\colon X\to \R$ a $\mu$-measurable function. Then for every $\e>0$ there exists a closed precompact set $K\subseteq X$ with $\mu(K)>1-\e$ and such that $f\vert_K$ is continuous.
\label{th:luzin}
\end{thrm}

As we could not find an exact reference to this specific version, we are including the proof.

\begin{thrm}
Every Borel probability measure, $\mu$, on a separable metric space $\Omega$ (not necessarily complete) satisfies the following regularity condition. Let $A$ be a $\mu$-measurable subset of $\Omega$. For every $\ve>0$ there exist a closed subset, $F$, and an open subset, $U$, of $\Omega$ such that
\[F\subseteq A\subseteq U\]
and
\[\mu(U\setminus F)<\ve.\]
\label{t:regularidadeFU}
\end{thrm}

\begin{proof}
Denote $\mathcal A$ the family of all Borel subsets of $\Omega$ satisfying the conclusion of the theorem: given $\ve>0$, there exist a closed set, $F$, and an open set, $U$, of $\Omega$ satisfying $F\subseteq A\subseteq U$ and $\mu(U\setminus F)<\ve$. It is easy to see that $\mathcal A$ forms a sigma-algebra which contains all closed subsets. Consequently, $\mathcal A$ contains all Borel sets. Since every $\mu$-measurable set differs from a suitable Borel set by a $\mu$-null set, we conclude. 
\end{proof}

\begin{proof}[Proof of Luzin's theorem \ref{th:luzin}]
Let $(\ve_n)$, $\ve_n>0$ be a summable sequence with $\sum_{n=1}^{\infty}\ve_n=\ve$. Enumerate the family of all open intervals with rational endpoints: $(a_n,b_n)$, $n\in\N$. For every $n$, use Th. \ref{t:regularidadeFU} to select closed sets $F_n\subseteq f^{-1}(a_n,b_n)$ and $F^{\prime}_n$, $F^{\prime}_n\cap f^{-1}(a_n,b_n)=\emptyset$ so that their union $\tilde F_n= F_n\cup F^\prime_n$ satisfies $\mu(\tilde F_n)>1-\ve_n/2$. 

The measure $\mu$ viewed as a Borel probability measure on the completion $\hat\Omega$ of the metric space $\Omega$ is regular, so there exists a compact set $Q\subseteq\hat\Omega$ with $\mu(Q)>1-\e/2$. 

The set
\[K=\bigcap_{n=1}^{\infty} \tilde F_n\cap Q\]
is closed and precompact in $\Omega$, and
satisfies $\mu(K)>1-\ve$. For each $n$, the set
$f^{-1}(a_n,b_n)\cap K$ is relatively open in $K$ because its compliment,
\[K\setminus f^{-1}(a_n,b_n) = F^{\prime}_n\cap K,\]
is closed.
\end{proof}

This simple proof is borrowed from \cite{LT}.

\section*{Concluding remarks}

The following question remains open. Let $\Omega$ be a separable complete metric space in which the $k$-NN classifier is universally consistent. Does it follow that $\Omega$ is sigma-finite dimensional in the sense of Nagata? 

A positive answer would imply, modulo the results of C\'erou and Guyader \cite{CG} and of Preiss \cite{preiss83}, that a separable metric space $\Omega$ satisfies the weak Lebesgue--Besicovitch differentiation property for every Borel sigma-finite locally finite measure if and only if $\Omega$ satisfies the strong Lebesgue--Besicovitch differentiation property for every Borel sigma-finite locally finite measure, which would answer an old question asked by Preiss in \cite{preiss83}. 

Most of this investigation appears as a part of the Ph.D. thesis of one of the authors \cite{kumari}.

The authors are thankful to the two anonymous referees of the paper, whose comments have helped to considerably improve the presentation. The remaining mistakes are of course authors' own.

\end{document}